\documentclass[10pt,reqno]{amsart}
\usepackage[hypertex]{hyperref}

\usepackage{amsmath,amsthm}
\usepackage{amssymb}
\usepackage{euscript}
\usepackage[all,tips]{xy}

\numberwithin{equation}{subsection}

\newcommand{\sqsp}{\renewcommand{\baselinestretch}{1.15}\tiny\normalsize}

\newcommand{\relphantom[1]}{\mathrel{\phantom{#1}}}
\newtheorem{theorem}[subsection]{Theorem}
\newtheorem{lemma}[subsection]{Lemma}

\newtheorem{corollary}[subsection]{Corollary}

\theoremstyle{definition}
\newtheorem{definition}[subsection]{Definition}
\newtheorem{example}[subsection]{Example}

\newcommand{\bk}{\mathbf{k}}
\newcommand{\bZ}{\mathbf{Z}}
\newcommand{\oct}{\mathbf{O}}

\newcommand{\xbar}{\overline{x}}
\newcommand{\ybar}{\overline{y}}
\newcommand{\zbar}{\overline{z}}

\newcommand{\muop}{\mu^{op}}

\newcommand{\bracket}{[,\ldots,]}
\newcommand{\form}{\langle,\rangle}

\DeclareMathOperator{\Hom}{Hom}
\begin{document}

\title{On $n$-ary Hom-Nambu and Hom-Nambu-Lie algebras}
\author{Donald Yau}

\begin{abstract}
It is observed that the category of $n$-ary Hom-Nambu(-Lie) algebras is closed under twisting by self-weak morphisms.  Constructions of ternary Hom-Nambu algebras from Hom-associative algebras, Hom-Lie algebras, ternary totally Hom-associative algebras, and Hom-Jordan triple systems are given.  Every multiplicative $n$-ary Hom-Nambu algebra gives rise to a sequence of Hom-Nambu algebras of exponentially higher arities.  Under some conditions, an $n$-ary Hom-Nambu(-Lie) algebra gives rise to an $(n-1)$-ary Hom-Nambu(-Lie) algebra.
\end{abstract}

\keywords{Hom-Nambu algebra, Hom-Nambu-Lie algebra, ternary totally Hom-associative algebra, Hom-Jordan triple system, Hom-Lie triple system, Hom-Lie algebra, Hom-associative algebra.}

\subjclass[2000]{17A40, 17A42, 17C50, 17D05, 17D10}

\address{Department of Mathematics\\
    The Ohio State University at Newark\\
    1179 University Drive\\
    Newark, OH 43055, USA}
\email{dyau@math.ohio-state.edu}

\date{\today}
\maketitle

\sqsp

%%%%%%%%%%%%%%%%%%%%%%
\section{Introduction}
%%%%%%%%%%%%%%%%%%%%%%

% n-ary algebras in Lie and Jordan algebras, geometry, math physics, M-branes

Algebras with $n$-ary compositions play important roles in Lie and Jordan theories, geometry, analysis, and physics.  For instance, Jordan triple systems \cite{jacobson1,lister,neher} give rise to $3$-graded Lie algebras through the TKK construction \cite{kantor1,koecher,tits}, from which most simple Lie algebras can be obtained.  Jordan triple systems also give rise to Lie triple systems through the Meyberg construction \cite{meyberg,meyberg2} (see \eqref{hjtshlts}).  On the other hand, Lie triple systems give rise to $\bZ/2\bZ$-graded Lie algebras \cite{jacobson1,lister}, which are exactly the kind of Lie algebras associated to symmetric spaces.  In geometry and analysis, various types of Jordan triple systems are used in the classifications of different classes of symmetric spaces \cite{bertram,chu,kaup1,kaup2,loos1,loos2}.

In physics, algebras with $n$-ary compositions appear in many different contexts.  Let us mention a few of them.  Ternary algebras can be used to construct solutions of the Yang-Baxter equation \cite{okubo} (section 8.5), which first appeared in statistical mechanics \cite{baxter,baxter2,yang}.  Nambu mechanics \cite{nambu} involves an $n$-ary product that satisfies the $n$-ary Nambu identity, which is an $n$-ary generalization of the Jacobi identity (see Definition \ref{def:homnambulie}).  Bagger-Lambert algebras \cite{bl} are ternary Nambu algebras with some extra structures, and they appear in the study of string theory and $M$-branes.  Ternary algebras are used in \cite{gun,gh,gh2} to construct superconformal algebras.

% n-ary Hom-algebras of Lie and associative types
% origin in Hom-Lie, other recent work

Generalizations of $n$-ary Nambu and Nambu-Lie algebras, called $n$-ary Hom-Nambu and Hom-Nambu-Lie algebras, were introduced in \cite{ams} by Ataguema, Makhlouf, and Silvestrov.  In these Hom-type algebras, the $n$-ary Nambu identity is relaxed using $n-1$ linear maps, called the twisting maps, resulting in the $n$-ary Hom-Nambu identity (Definition \ref{def:homnambulie}).  When these twisting maps are all equal to the identity map, one recovers $n$-ary Nambu and Nambu-Lie algebras.

% Hom-algebras

Binary Hom-Lie algebras, or just Hom-Lie algebras, were introduced in \cite{hls} to describe the structures on some $q$-deformations of the Witt and the Virasoro algebras.  The associative version of a Hom-Lie algebra, called a Hom-associative algebra, was defined in \cite{ms}.  They are related in roughly the same way in which Lie and associative algebras are related \cite{ms,yau}.  These Hom-type algebras are defined by twisting the defining identities (the Jacobi identity and associativity) by a linear twisting map.  Other Hom-type generalizations of familiar algebras can be found in \cite{mak} - \cite{ms4} and \cite{yau3,yau4,yau12}.  Hom-type analogues of quantum groups and the Yang-Baxter equations have been studied in \cite{yau5} - \cite{yau11}.

% constructions of ternary Hom-Nambu-Lie

Several general recipes for constructing $n$-ary Hom-type algebras are known.  A method for constructing $n$-ary Hom-Nambu(-Lie) algebras from $n$-ary Nambu(-Lie) algebras and endomorphisms was given in \cite{ams}.  It is a generalization of a result due to the author \cite{yau2}.  More recently, in \cite{ams1} it is shown that certain $q$-deformations of the ternary Virasoro-Witt algebra \cite{cfz} are ternary Hom-Nambu-Lie algebras.  Using a ternary product first used in \cite{almy}, it is shown in \cite{ams2} that one can obtain a ternary Hom-Nambu-Lie algebra from a Hom-Lie algebra with a compatible linear map and a trace function.

% purpose: how different categories of n-ary Hom-Nambu(-Lie) are related

The purpose of this paper is to study the relationships between Hom-Nambu(-Lie) algebras of different arities and between Hom-Nambu algebras and other Hom-type algebras.  A description of the rest of this paper follows.

In section \ref{sec:twist} we recall the definitions of $n$-ary Hom-Nambu and Hom-Nambu-Lie algebras.  It is shown that the category of $n$-ary Hom-Nambu(-Lie) algebras is closed under twisting by self-weak morphisms (Theorem \ref{thm:twist}).  The same kind of closure property holds in other categories of Hom-type algebras.  We stress that this is a unique property for Hom-type algebras, since ordinary algebras are usually not closed under such twistings.  Restricting such twistings to $n$-ary Nambu or Nambu-Lie algebras, one recovers the twisting result (Corollary \ref{cor3:twist}) in \cite{ams}.

In section \ref{sec:homtriple} we study ternary Hom-algebras of associative, Jordan, and Lie types.  We recall the notion of a ternary totally Hom-associative algebra from \cite{ams} and define Hom-Jordan and Hom-Lie triple systems, which generalize Jordan and Lie triple systems.  Hom-Lie triple systems are ternary Hom-Nambu algebras whose triple product is left anti-symmetric and satisfies the ternary Jacobi identity (see \eqref{homltsid}).  Each of these categories of ternary Hom-algebras is closed under twisting by self-weak morphisms (Theorem \ref{thm:twisttriple}).  Using a special case of this twisting property (Corollary \ref{cor2:twisttriple}), we discuss how Hom-Jordan triple systems arise from Jordan algebras (Corollary \ref{cor:jhjts}) and alternative algebras (Corollary \ref{cor:althjts}).  We also discuss how Hom-Lie triple systems arise from Maltsev algebras (Corollary \ref{cor:maltsevhlts}) and alternative algebras (Corollary \ref{cor:althlts}).

In section \ref{sec:meyberg} we show that Hom-Lie triple systems, and hence ternary Hom-Nambu algebras, arise from Hom-Jordan triple systems, ternary totally Hom-associative algebras, Hom-associative algebras, and Hom-Lie algebras.  In particular, it is shown that every ternary totally Hom-associative algebra with equal twisting maps has an underlying Hom-Jordan triple system (Theorem \ref{thm:atsjts}).  Generalizing the Meyberg construction, it is shown that every Hom-Jordan triple system with equal twisting maps has an underlying Hom-Lie triple system (Theorem \ref{thm:jtslts}).  Combining these two results, it follows that every ternary totally Hom-associative algebra with equal twisting maps has an underlying Hom-Lie triple system (Corollary \ref{cor:atslts}).  Then we shown that Hom-associative and multiplicative Hom-Lie algebras give rise to ternary totally Hom-associative algebras and multiplicative Hom-Lie triple systems (Theorems \ref{thm:haats} and \ref{thm:hllts}).  Combining Corollary \ref{cor:atslts} and Theorem \ref{thm:haats}, it follows that every Hom-associative algebra has an underlying Hom-Lie triple system (Corollary \ref{cor:hahnambu}).

In section \ref{sec:arity} it is shown that multiplicative $n$-ary Hom-Nambu algebras give rise to Hom-Nambu algebras of exponentially higher arities.  In particular, it is shown that every multiplicative $n$-ary Hom-Nambu algebra $L$ gives rise to a multiplicative $(2n-1)$-ary Hom-Nambu algebra $L'$.  The $(2n-1)$-ary product in $L'$  involves a two-fold composition of the $n$-ary product in $L$ and $n-1$ copies of the twisting map (Theorem \ref{thm:higher}).  This observation is inspired by Jacobson's original definition of a Lie triple system \cite{jacobson1}, which was defined there as a submodule of an associative algebra that is closed under the iterated commutator bracket $[[x,y],z]$.  Iterating Theorem \ref{thm:higher}, we obtain from a multiplicative $n$-ary Hom-Nambu algebra $L$ a sequence of Hom-Nambu algebras $L^k$ of arities $2^k(n-1) + 1$ for $k \geq 0$ (Corollary \ref{cor1:higher}).

In section \ref{sec:lower} it is shown that, under some conditions, an $n$-ary (with $n \geq 3$) Hom-Nambu(-Lie) algebra yields a reduced Hom-Nambu(-Lie) algebra of arity $n-1$ (Theorem \ref{thm:lower}).  In particular, for an $n$-ary Hom-Nambu-Lie algebra $L$, the product in the reduced $(n-1)$-ary Hom-Nambu-Lie algebra takes the form $[a,\ldots]$, where $a$ is a fixed element of the first twisting map $\alpha_1$ (Corollary \ref{cor1:lower}).  Iterating these results, we obtain a way to reduce an $n$-ary Hom-Nambu(-Lie) algebra to an $(n-k)$-ary Hom-Nambu(-Lie) algebra for $k \in \{1,\ldots,n-2\}$ (Corollary \ref{cor2:lower} and Corollary \ref{cor3:lower}).  Combining Corollary \ref{cor1:lower} with a result in \cite{ams2}, a non-trivial method for replacing the structure maps in a Hom-Lie algebra is given (Corollary \ref{cor5:lower}).

%%%%%%%%%%%%%%%%%%%%%%%%%%%%%%%%%%%%%%%%%%%
\section{Twisting Hom-Nambu(-Lie) algebras}
\label{sec:twist}
%%%%%%%%%%%%%%%%%%%%%%%%%%%%%%%%%%%%%%%%%%%

The purpose of this section is to observe that the category of $n$-ary Hom-Nambu(-Lie) algebras is closed under twisting by self-weak morphisms.  Some consequences and examples are also discussed.

%%%%%%%%%%%%%%%%%%%%%%%%
\subsection{Conventions}

Throughout this paper we work over a fixed field $\bk$ of characteristic $0$.  If $V$ is a $\bk$-module and $f \colon V \to V$ is a linear map, then $f^n$ denotes the composition of $n$ copies of $f$ with $f^0 = Id$.  For elements $x_i,\ldots,x_j \in V$, we use the abbreviations
\begin{equation}
\label{xij}
\begin{split}
x_{i,j} &= (x_i,x_{i+1},\ldots,x_j),\\
f(x_{i,j}) &= (f(x_i),\ldots, f(x_j))\\
\end{split}
\end{equation}
when $i \leq j$; for $i > j$, $x_{i,j}$ and $f(x_{i,j})$ denote the empty sequence.

Let us begin with the following basic definitions.

% definitions
\begin{definition}
\label{def:nhomalgebra}
Let $n \geq 2$ be an integer.
\begin{enumerate}
\item
An \textbf{$n$-ary Hom-algebra} $(V,\bracket,\alpha=(\alpha_1,\ldots,\alpha_{n-1}))$ consists of a $\bk$-module $V$, an $n$-linear map $\bracket \colon V^{\otimes n} \to V$, and linear maps $\alpha_i \colon V \to V$ for $i = 1, \ldots , n-1$, called the \textbf{twisting maps}.
\item
An $n$-ary Hom-algebra $(V,\bracket,\alpha)$ is said to be \textbf{multiplicative} if the twisting maps are all equal, i.e., $\alpha_1 = \cdots = \alpha_{n-1} = \alpha$, and $\alpha \circ \bracket = \bracket \circ \alpha^{\otimes n}$.
\item
A \textbf{weak morphism} $f \colon V \to U$ of $n$-ary Hom-algebras is a linear map of the underlying $\bk$-modules such that $f \circ \bracket_V = \bracket_U \circ f^{\otimes n}$.  A \textbf{morphism} of $n$-ary Hom-algebras is a weak morphism such that $f \circ (\alpha_i)_V = (\alpha_i)_U \circ f$ for $i = 1, \ldots n-1$ \cite{ams}.
\item
The $n$-ary product $\bracket$ in an $n$-ary Hom-algebra $V$ is said to be \textbf{anti-symmetric} if
\[
[x_1,\ldots,x_n] = \epsilon(\sigma)[x_{\sigma(1)},\ldots,x_{\sigma(n)}]
\]
for all $x_i \in V$ and permutations $\sigma$ on $n$ letters, where $\epsilon(\sigma)$ is the signature of $\sigma$.
\end{enumerate}
\end{definition}

When all the twisting maps are equal in an $n$-ary Hom-algebra $V$, as in the multiplicative case, we will write it as $(V,\bracket,\alpha)$, where $\alpha$ is the common value of the twisting maps.

As is customary, an $n$-ary Hom-algebra is called \textbf{binary} or \textbf{ternary} if $n=2$ or $3$, respectively.  An $n$-ary algebra in the usual sense is a $\bk$-module $V$ with an $n$-linear map $\bracket \colon V^{\otimes n} \to V$.  We consider an $n$-ary algebra $(V,\bracket)$ also as an $n$-ary Hom-algebra $(V,\bracket,Id)$ in which all $n-1$ twisting maps are the identity map.  Also, in this case a weak morphism is the same thing as a morphism, which agrees with the usual definition of a morphism of $n$-ary algebras.

Next we recall the definition of an $n$-ary Hom-Nambu(-Lie) algebra from \cite{ams}.

\begin{definition}
\label{def:homnambulie}
Let $(V,\bracket,\alpha=(\alpha_1,\ldots,\alpha_{n-1}))$ be an $n$-ary Hom-algebra.
\begin{enumerate}
\item
The \textbf{$n$-ary Hom-Jacobian} of $V$ is the $(2n-1)$-linear map $J^n_V \colon V^{\otimes 2n-1} \to V$ defined as (using the shorthand in \eqref{xij})
\begin{equation}
\label{homjacobian}
\begin{split}
J^n_V(x_{1,n-1};y_{1,n})
&= [\alpha_1(x_1),\ldots,\alpha_{n-1}(x_{n-1}),[y_{1,n}]]\\
&\relphantom{} - \sum_{i=1}^n [\alpha_1(y_1),\ldots,\alpha_{i-1}(y_{i-1}), [x_{1,n-1},y_i], \alpha_i(y_{i+1}),\ldots,\alpha_{n-1}(y_n)]
\end{split}
\end{equation}
for $x_1, \ldots, x_{n-1}, y_1, \ldots , y_n \in V$.
\item
An \textbf{$n$-ary Hom-Nambu algebra} is an $n$-ary Hom-algebra $V$ that satisfies the \textbf{$n$-ary Hom-Nambu identity} $J^n_V = 0$.
\item
An \textbf{$n$-ary Hom-Nambu-Lie algebra} is an $n$-ary Hom-Nambu algebra in which the $n$-ary product $\bracket$ is anti-symmetric.
\item
Multiplicativity and (weak) morphisms for $n$-ary Hom-Nambu(-Lie) algebras are defined using the underlying $n$-ary Hom-algebras.
\end{enumerate}
\end{definition}

When the twisting maps are all equal to the identity map, $n$-ary Hom-Nambu algebras and $n$-ary Hom-Nambu-Lie algebras are the usual \textbf{$n$-ary Nambu algebras} and \textbf{$n$-ary Nambu-Lie algebras} \cite{filippov,nambu}.  In this case, the $n$-ary Hom-Jacobian $J^n_V$ is called the \textbf{$n$-ary Jacobian} \cite{filippov}, and the $n$-ary Hom-Nambu identity $J^n_V = 0$ is called the \textbf{$n$-ary Nambu identity}.

There is a way to construct $n$-ary Hom-Nambu(-Lie) algebras from $n$-ary Nambu(-Lie) algebras by twisting along morphisms \cite{ams} (Theorem 3.4).  That result is a variation of a result by the author \cite{yau2} (Theorem 2.3).  We will establish a more general twisting result below and obtain \cite{ams} (Theorem 3.4) as a special case.  To this end, we need the following observations.

\begin{lemma}
\label{lem:twist}
Let $(V,\bracket,\alpha)$ be an $n$-ary Hom-algebra and $\beta \colon V \to V$ be a weak morphism.  Define the $n$-ary product $\bracket_\beta = \beta \circ \bracket$ and the twisting maps $\beta\alpha = (\beta\alpha_1,\ldots,\beta\alpha_{n-1})$.  Consider the $n$-ary Hom-algebra
\begin{equation}
\label{vbeta}
V_\beta = (V,\bracket_\beta,\beta\alpha).
\end{equation}
Then the following statements hold.
\begin{enumerate}
\item
$J^n_{V_\beta} = \beta^2 \circ J^n_V$.
\item
If $\bracket$ is anti-symmetric, then so is $\bracket_\beta$.
\item
If $V$ is multiplicative and $\beta\alpha = \alpha\beta$, then $V_\beta$ is multiplicative.
\end{enumerate}
\end{lemma}

\begin{proof}
All three statements are immediate from the definitions.
\end{proof}

Using Lemma \ref{lem:twist}, we obtain immediately the following twisting result.

\begin{theorem}
\label{thm:twist}
If $(V,\bracket,\alpha)$ is an $n$-ary Hom-Nambu(-Lie) algebra and $\beta \colon V \to V$ is a weak morphism, then $V_\beta$ in \eqref{vbeta} is also an $n$-ary Hom-Nambu(-Lie) algebra.  Moreover, if $V$ is multiplicative and $\beta\alpha = \alpha\beta$, then $V_\beta$ is also multiplicative.
\end{theorem}

Note that Theorem \ref{thm:twist}, as well as Corollaries \ref{cor1:twist} - \ref{cor3:twist} below, have obvious analogues for $n$-ary Hom-Lie algebras and $n$-ary totally (or partially) Hom-associative algebras, as defined in \cite{ams}.

Theorem \ref{thm:twist} says that the category of $n$-ary Hom-Nambu(-Lie) algebras is closed under twisting by self-weak morphisms.  The same kind of closure under twisting result can be established for most other kinds of Hom-type algebras. This property is unique to Hom-type algebras, as the usual notions of algebras are usually not closed under twisting by morphisms.

We record the following special cases of Theorem \ref{thm:twist}.  First, if $V$ is multiplicative and $\beta = \alpha$, then we have the following twisting result.

\begin{corollary}
\label{cor1:twist}
If $(V,\bracket,\alpha)$ is a multiplicative $n$-ary Hom-Nambu(-Lie) algebra, then so is $V_\alpha = (V,\bracket_\alpha = \alpha\circ\bracket,\alpha^2)$.
\end{corollary}

If we use Corollary \ref{cor1:twist} repeatedly, then we obtain a sequence of derived Hom-Nambu(-Lie) algebras as follows.

\begin{corollary}
\label{cor2:twist}
If $(V,\bracket,\alpha)$ is a multiplicative $n$-ary Hom-Nambu(-Lie) algebra, then so is $V_k = (V,\bracket_k = \alpha^{2^k-1}\circ\bracket,\alpha^{2^k})$ for each $k \geq 0$.
\end{corollary}

On the other hand, if the twisting maps $\alpha_i$ of $V$ are all equal to the identity map in Theorem \ref{thm:twist}, then we recover the following twisting result from \cite{ams} (Theorem 3.4).

\begin{corollary}
\label{cor3:twist}
Let $(V,\bracket)$ be an $n$-ary Nambu(-Lie) algebra and $\beta \colon V \to V$ be a morphism.  Then $V_\beta = (V,\bracket_\beta = \beta \circ \bracket,\beta)$ is a multiplicative $n$-ary Hom-Nambu(-Lie) algebra.
\end{corollary}

Some examples of ternary Hom-Nambu(-Lie) algebras constructed using Corollary \ref{cor3:twist} can be found in \cite{ams1,ams2,ams}.  In the rest of this section, we give some examples of ternary Hom-Nambu algebras.

% examples: JTS, ATS, symmetric bilinear form (Okubo's p.113)

%%%%%%%%%%%%%%%%%%%
\begin{example}
\label{ex:okubo}
Let $V$ be a $\bk$-module and $\form \colon V^{\otimes 2} \to \bk$ be a symmetric bilinear form.  Then for any scalar $\lambda \in \bk$, the triple product
\begin{equation}
\label{formproduct}
[xyz] = \lambda(\langle y,z\rangle x - \langle z,x\rangle y)
\end{equation}
gives a ternary Nambu algebra $(V,[,,])$ \cite{okubo} (p.113).  In fact, it is a Lie triple system \cite{lister}, which is always a ternary Nambu algebra.  See Definition \ref{def:hjts}.

Suppose $\alpha \colon V \to V$ is a linear map that is invariant with respect to the bilinear form $\form$ in the sense that
\begin{equation}
\label{invariant}
\langle x,y\rangle = \langle \alpha(x),\alpha(y)\rangle
\end{equation}
for all $x,y \in V$.  Then $\alpha$ is also a morphism on the ternary Nambu algebra $(V,[,,])$.  By Corollary \ref{cor3:twist} there is a multiplicative ternary Hom-Nambu algebra $V_\alpha = (V,[,,]_\alpha=\alpha\circ[,,],\alpha)$ with
\begin{equation}
\label{formalphaproduct}
[xyz]_\alpha = \lambda(\langle y,z\rangle \alpha(x) - \langle z,x\rangle \alpha(y)).
\end{equation}
Note that $V_\alpha$ is usually not a ternary Hom-Nambu-Lie algebra because the triple product $[,,]_\alpha$ is in general not anti-symmetric.  Moreover, $V_\alpha$ is usually not a ternary Nambu algebra, as the next example illustrates.
\qed
\end{example}
%%%%%%%%%%%%%%%%%%%

%%%%%%%%%%%%%%%%%%%
\begin{example}
\label{ex:fermion}
For $N \geq 2$ consider the $2N$-dimensional vector space $V$ spanned by the fermionic annihilation operators $a_{-j}$ and the creation operators $a_{+j}$ for $1 \leq j \leq N$ with the symmetric bilinear form determined by
\begin{equation}
\label{fermion}
\langle a_{\mp j},a_{\pm k}\rangle = \delta_{jk},\quad
\langle a_{\pm j}, a_{\pm k}\rangle = 0.
\end{equation}
As in Example \ref{ex:okubo}, we have a ternary Nambu algebra (in fact, a Lie triple system) $(V,[,,])$ with the triple product in \eqref{formproduct} \cite{okubo} (p.113).

For each $j \in \{1,\ldots,N\}$, let $\eta_j$ be a non-zero scalar in $\bk$.  Let $\alpha \colon V \to V$ be the linear automorphism defined by
\begin{equation}
\label{alphaf}
\alpha(a_{\pm j}) = \eta_j^{\pm 1}a_{\pm j}.
\end{equation}
Then $\alpha$ is invariant with respect to the symmetric bilinear form in \eqref{fermion} in the sense of \eqref{invariant}.  By Corollary \ref{cor3:twist} there is a multiplicative ternary Hom-Nambu algebra $V_\alpha = (V,[,,]_\alpha,\alpha)$, where $[,,]_\alpha = \alpha\circ[,,]$ is the triple product in \eqref{formalphaproduct}.  In fact, this triple product is completely determined by
\begin{equation}
\label{fsystem}
\begin{split}
[a_{\pm i}, a_{\pm j}, a_{\mp k}]_\alpha &= \lambda(\delta_{jk}\eta_i^{\pm 1}a_{\pm i} - \delta_{ki}\eta_j^{\pm 1} a_{\pm j}),\\
[a_{\pm i}, a_{\mp j}, a_{\pm k}]_\alpha &= \lambda\delta_{jk}\eta_i^{\pm 1} a_{\pm i},\\
[a_{\mp i}, a_{\pm j}, a_{\pm k}]_\alpha &= - \lambda\delta_{ki}\eta_j^{\pm 1} a_{\pm j}.
\end{split}
\end{equation}
Moreover, $(V,[,,]_\alpha)$ is not a ternary Nambu algebra because the ternary Jacobian for $(V,[,,]_\alpha)$ gives
\[
J^3(a_{+1}a_{+2};a_{-2}a_{+2}a_{-2}) = \lambda^2(\eta_1\eta_2^{-1} - \eta_1^2)a_{+1}.
\]
So the ternary Nambu identity is not satisfied as long as $\eta_1 \not= \eta_2^{-1}$ and $\lambda \not= 0$.

Finally, note that the middle relation in \eqref{fsystem} (when $\lambda = 2$ and $\eta_i = 1$) is a relation in the para-fermionic system \cite{ok}.
\qed
\end{example}
%%%%%%%%%%%%%%%%%%%

%%%%%%%%%%%%%%%%%%%%%%%%%%%%%%%%%%%%
\section{Ternary Hom-algebras of associative, Jordan, and Lie types}
\label{sec:homtriple}
%%%%%%%%%%%%%%%%%%%%%%%%%%%%%%%%%%%%

In this section, we discuss the three types of ternary Hom-algebras in the section title.  They generalize the usual triple systems of associative, Jordan, and Lie types.  We observe that each of these types of ternary Hom-algebras is closed under twisting by self-weak morphisms (Theorem \ref{thm:twisttriple}) and discuss some examples using a special case of this twisting result.  Our main interest in these ternary Hom-algebras lies in the fact that they give rise to ternary Hom-Nambu algebras, as we discuss in the next section.

Let us first recall the following ternary generalization of Hom-associative algebras from \cite{ams}

\begin{definition}
\label{def:homats}
A \textbf{ternary totally Hom-associative algebra} \cite{ams} is a ternary Hom-algebra $(A,(,,),\alpha=(\alpha_1,\alpha_2))$ that satisfies \textbf{ternary Hom-associativity}
\begin{equation}
\label{thomass}
((uvw)\alpha_1(x)\alpha_2(y)) = (\alpha_1(u)(vwx)\alpha_2(y)) = (\alpha_1(u)\alpha_2(v)(wxy))
\end{equation}
for all $u,v,w,x,y \in A$.
\end{definition}

When both twisting maps $\alpha_i$ are equal to the identity map, a ternary totally Hom-associative algebra becomes a \textbf{ternary ring} in the sense of \cite{lister2}.  In this case, the condition \eqref{thomass} is called \textbf{ternary associativity}.  Ternary rings are also known as associative triple systems in the literature.

Next we define the Jordan and Lie analogues of ternary totally Hom-associative algebras.

% (multiplicative) Hom-ATS, JTS, LTS
\begin{definition}
\label{def:hjts}
\begin{enumerate}
\item
A \textbf{Hom-triple system} is a ternary Hom-algebra $(V,\{,,\},\alpha=(\alpha_1,\alpha_2))$.  Multiplicativity, weak morphisms, and morphisms are defined as in Definition \ref{def:nhomalgebra}.
\item
A \textbf{Hom-Jordan triple system} is a Hom-triple system $(J,\{,,\},\alpha)$ that satisfies
\begin{equation}
\label{outersymmetry}
\{xyz\} = \{zyx\} \quad\text{(outer-symmetry)}
\end{equation}
and the \textbf{Hom-Jordan triple identity}
\begin{equation}
\label{homjtsid}
\begin{split}
\{\alpha_1(x)\alpha_2(y)\{uvw\}\} &- \{\alpha_1(u)\alpha_2(v)\{xyw\}\}\\
&= \{\{xyu\}\alpha_1(v)\alpha_2(w)\} - \{\alpha_1(u)\{yxv\}\alpha_2(w)\}
\end{split}
\end{equation}
for all $u,v,w,x,y \in J$.
\item
A \textbf{Hom-Lie triple system} is a Hom-triple system $(L,[,,],\alpha)$ that satisfies
\begin{equation}
\label{homltsid}
\begin{split}
[uvw] &= -[vuw] \quad\text{(left anti-symmetry)},\\
0 &= [uvw] + [wuv] + [vwu] \quad\text{(ternary Jacobi identity)},
\end{split}
\end{equation}
and the ternary Hom-Nambu identity $J^3_L = 0$, i.e.,
\begin{equation}
\label{homnambu}
[\alpha_1(x)\alpha_2(y)[uvw]] = [[xyu]\alpha_1(v)\alpha_2(w)] + [\alpha_1(u)[xyv]\alpha_2(w)] + [\alpha_1(u)\alpha_2(v)[xyw]]
\end{equation}
for all $u,v,w,x,y \in J$.
\end{enumerate}
\end{definition}

% remark: Hom-ATS as ternary totally hom-associative algebra
% Hom-LTS as ternary Hom-Nambu (but not Hom-Nambu-Lie) algebra

When the twisting maps $\alpha_i$ are both equal to the identity map, we recover the usual notions of a \textbf{Jordan triple system} \cite{meyberg} and a \textbf{Lie triple system} \cite{lister}.  So Jordan and Lie triple systems are examples of multiplicative Hom-Jordan and Hom-Lie triple systems, respectively.  In the special case when both twisting maps are the identity map, we call \eqref{homjtsid} the \textbf{Jordan triple identity}.

A Hom-Lie triple system is automatically a ternary Hom-Nambu algebra as in Definition \ref{def:homnambulie}.  Note that a Hom-Lie triple system, whose triple product is only assumed to be left antisymmetric, is in general not a ternary Hom-Nambu-Lie algebra.  In particular, all the constructions of Hom-Lie triple systems in this paper give examples of ternary Hom-Nambu algebras that are usually not ternary Hom-Nambu-Lie algebras.  The reader is advised not to confuse a Hom-Lie triple system with a ternary Hom-Lie algebra as defined in \cite{ams}.

The following result is the analogue of Theorem \ref{thm:twist} for the ternary Hom-algebras above.  It is proved in essentially the same way, i.e., by applying either $\beta$ or $\beta^2$ to the defining identities in Definitions \ref{def:homats} and \ref{def:hjts}.

% twisting results

\begin{theorem}
\label{thm:twisttriple}
Let $(V,(,,),\alpha=(\alpha_1,\alpha_2))$ be a ternary totally Hom-associative algebra (resp., Hom-Jordan or Hom-Lie triple system) and $\beta \colon V \to V$ be a weak morphism.  Then the ternary Hom-algebra
\[
V_\beta = (V,(,,)_\beta = \beta \circ (,,),\beta\alpha = (\beta\alpha_1,\beta\alpha_2))
\]
is a ternary totally Hom-associative algebra (resp., Hom-Jordan or Hom-Lie triple system).  If, in addition, $V$ is multiplicative and $\beta\alpha = \alpha\beta$, then $V_\beta$ is multiplicative.
\end{theorem}

The following result is the special case of Theorem \ref{thm:twisttriple} when $\beta = \alpha$.

% multiplicative hom-triple twisted along its own twisting map
\begin{corollary}
\label{cor1:twisttriple}
Let $(V,(,,),\alpha)$ be a multiplicative ternary totally Hom-associative algebra (resp., Hom-Jordan or Hom-Lie triple system).  Then $V_\alpha = (V,(,,)_\alpha = \alpha \circ (,,),\alpha^2)$ is also  a multiplicative ternary totally Hom-associative algebra (resp., Hom-Jordan or Hom-Lie triple system).
\end{corollary}

Using Corollary \ref{cor1:twisttriple} repeatedly, we obtain a sequence of derived ternary Hom-algebras.

\begin{corollary}
\label{cor1.5:twisttriple}
Let $(V,(,,),\alpha)$ be a multiplicative ternary totally Hom-associative algebra (resp., Hom-Jordan or Hom-Lie triple system).  Then for each $k \geq 0$,
\[
V_k = (V,\alpha^{2^k-1} \circ (,,),\alpha^{2^k})
\]
is also a multiplicative ternary totally Hom-associative algebra (resp., Hom-Jordan or Hom-Lie triple system).
\end{corollary}

This following result is the special case of Theorem \ref{thm:twisttriple} when the twisting maps $\alpha_i$ are both equal to the identity map on $V$.  The case concerning a ternary ring is also a special case of \cite{ams} (Theorem 3.6).

% twisting A/J/LTS into Hom-A/J/LTS
\begin{corollary}
\label{cor2:twisttriple}
Let $(V,(,,))$ be a ternary ring (resp., Jordan or Lie triple system) and $\beta \colon V \to V$ be a morphism of ternary algebras.  Then
$V_\beta = (V,(,,)_\beta = \beta \circ (,,),\beta)$ is a multiplicative ternary totally Hom-associative algebra (resp., Hom-Jordan or Hom-Lie triple system).
\end{corollary}

% examples
The rest of this section contains consequences of Corollary \ref{cor2:twisttriple} and examples of ternary totally Hom-associative algebras, Hom-Jordan triple systems, and Hom-Lie triple systems.  Some examples of ternary totally (or partially) Hom-associative algebras can be found in \cite{ams}.  We begin with examples of ternary totally Hom-associative algebras, for which we use the following observation.

\begin{corollary}
\label{cor:homats}
Let $(A,\mu)$ be an associative algebra and $B$ be a submodule of $A$ that is closed under the triple product $(,,) = \mu \circ (\mu \otimes Id)$.  Then $(B,(,,))$ is a ternary ring.  Moreover, if $B$ is a sub-algebra of $A$ and $\alpha \colon B \to B$ is an algebra morphism (i.e., $\alpha\circ\mu|_B = \mu|_B \circ \alpha^{\otimes 2}$), then $B_\alpha = (B,(,,)_\alpha,\alpha)$ is a multiplicative ternary totally Hom-associative algebra, where $(,,)_\alpha = \alpha \circ (,,)$.
\end{corollary}

\begin{proof}
It is clear that the triple product $(xyz) = (xy)z$ satisfies ternary associativity, so $(B,(,,))$ is a ternary ring.  The second assertion follows from Corollary \ref{cor2:twisttriple}, since $\alpha$ is multiplicative with respect to the triple product $(,,) = \mu \circ (\mu \otimes Id)$.
\end{proof}

%%%%%%%%%%%%%%%%%%%%%
\begin{example}
\label{ex1:homats}
Let $(A,\mu)$ be an associative algebra, $\theta \colon A \to A$ be an algebra morphism, and $B$ be the $(-1)$-eigenspace of $\theta$, i.e.,
\[
B = \{a \in A \colon \theta(a) = -a\}.
\]
Then $B$ is closed under the triple product $(,,) = \mu \circ (\mu \otimes Id)$, so $(B,(,,))$ is a ternary ring (Corollary \ref{cor:homats}).  If $\alpha \colon A \to A$ is another algebra morphism such that $\theta \circ \alpha = \alpha \circ \theta$, then $\alpha$ restricts to a map on $B$.  Therefore, by Corollary \ref{cor:homats} again, we have a multiplicative ternary totally Hom-associative algebra $B_\alpha = (B,(,,)_\alpha,\alpha)$ with $(,,)_\alpha = \alpha \circ \mu \circ (\mu \otimes Id)$.
\qed
\end{example}
%%%%%%%%%%%%%%%%%%%%%

%%%%%%%%%%%%%%%%%%%%%
\begin{example}
\label{ex2:homats}
Let $n>1$ be a positive integer and $A$ be the associative algebra consisting of polynomials over $\bk$ in $n$ variables $x_1, \ldots , x_n$ with $0$ constant term.  For the algebra morphism $\theta \colon A \to A$ defined by
\[
\theta(x_i) = -x_i \quad \text{for $1 \leq i \leq n$},
\]
the $(-1)$-eigenspace $B$ consists of the odd polynomials, i.e., linear combinations of monomials of odd total degrees.  By Example \ref{ex1:homats}, $B$ is a ternary ring.  This example (when $n=2$) is from \cite{lister2} p.47.

For each $i$, let $r_i > 1$ be a positive odd integer.  Then the algebra morphism $\alpha \colon A \to A$ defined by
\[
\alpha(x_i) = x_i^{r_i} \quad \text{for $1 \leq i \leq n$}
\]
commutes with $\theta$.  By Example \ref{ex1:homats}, $B_\alpha = (B,(,,)_\alpha,\alpha)$ with $(,,)_\alpha = \alpha \circ \mu \circ (\mu \otimes Id)$ is a multiplicative ternary totally Hom-associative algebra.  Finally, observe that $(B,(,,)_\alpha)$ is not a ternary ring.  In fact, we have
\[
((x_1x_1x_1)_\alpha x_2x_2)_\alpha = x_1^{3r_1^2}x_2^{2r_2} \not= x_1^{r_1(r_1+2)} x_2^{2r_2^2} = (x_1x_1(x_1x_2x_2)_\alpha)_\alpha,
\]
so $(B,(,,)_\alpha)$ does not satisfy ternary associativity.
\qed
\end{example}
%%%%%%%%%%%%%%%%%%%%%

%%%%%%%%%%%%%%%%%%%%%
\begin{example}
\label{ex3:homats}
Let $V$ and $W$ be two $\bk$-modules.  The direct sum
\[
A = \Hom(V,W) \oplus \Hom(W,V)
\]
becomes a ternary ring when equipped with the triple product
\[
(f_1 \oplus g_1, f_2 \oplus g_2, f_3 \oplus g_3) = (f_3g_2f_1) \oplus (g_3f_2g_1)
\]
for $f_i \in \Hom(V,W)$ and $g_j \in \Hom(W,V)$.  This example is from \cite{lister2} p.46.

Suppose $\beta \colon V \to V$ and $\gamma \colon W \to W$ are linear automorphisms.  Define the map $\alpha \colon A \to A$ by
\[
\alpha(f \oplus g) = (\gamma^{-1}f\beta) \oplus (\beta^{-1}g\gamma)
\]
for $f \in \Hom(V,W)$ and $g \in \Hom(W,V)$.  Then $\alpha$ is a morphism of ternary algebras.  In fact, $\alpha$ is an automorphism with inverse
\[
\alpha^{-1}(f \oplus g) = (\gamma f \beta^{-1}) \oplus (\beta g \gamma^{-1}).
\]
By Corollary \ref{cor2:twisttriple} we have a multiplicative ternary totally Hom-associative algebra $A_\alpha = (A,(,,)_\alpha,\alpha)$ with $(,,)_\alpha = \alpha \circ (,,)$.  Moreover, for $a_i \in A$ ($1 \leq i \leq 5$), one can check that
\[
((a_1a_2a_3)_\alpha a_4 a_5)_\alpha \not= (a_1a_2(a_3a_4a_5)_\alpha)_\alpha
\]
in general.  Therefore, $(A,(,,)_\alpha)$ is not a ternary ring.
\qed
\end{example}
%%%%%%%%%%%%%%%%%%%%%

The next few examples are about Hom-Jordan triple systems.

%%%%%%%%%%%%%%%%%%%%%
\begin{example}
\label{ex1:hjts}
Let $V$ be a $\bk$-module and $\langle,\rangle \colon V^{\otimes 2} \to \bk$ be a symmetric bilinear form.  Then $(V,\{,,\})$ is a Jordan triple system with triple product
\begin{equation}
\label{jhjproduct}
\{xyz\} = \lambda(\langle x,y\rangle z + \langle y,z\rangle x - \langle z,x\rangle y)
\end{equation}
for all $x,y,z \in V$, where $\lambda \in \bk$ is any fixed scalar.  Suppose $\alpha \colon V \to V$ is a linear map that is invariant with respect to the symmetric bilinear form in the sense that
\[
\langle \alpha(x),\alpha(y)\rangle = \langle x,y\rangle
\]
for all $x,y \in V$.  Then $\alpha \colon (V,\{,,\}) \to (V,\{,,\})$ is a morphism of Jordan triple systems.  By Corollary \ref{cor2:twisttriple} we have a multiplicative Hom-Jordan triple system $V_\alpha = (V,\{,,\}_\alpha,\alpha)$, where $\{,,\}_\alpha = \alpha\circ\{,,\}$.

If we use $\lambda/2$ rather than $\lambda$ in \eqref{jhjproduct}, then the triple product
\[
[xyz]_\alpha = \{xyz\}_\alpha - \{yxz\}_\alpha = \lambda(\langle y,z\rangle \alpha(x) - \langle z,x\rangle \alpha(y))
\]
allows us to recover the ternary Hom-Nambu product in \eqref{formalphaproduct}.  This is not a coincidence.  In fact, a Hom-Jordan triple system with equal twisting maps always gives rise to a Hom-Lie triple system, as we show in Theorem \ref{thm:jtslts} below.  Finally, since Example \ref{ex:fermion} is a special case of Example \ref{ex:okubo}, the ternary Hom-Nambu algebra in that Example is also induced by a Hom-Jordan triple system.
\qed
\end{example}
%%%%%%%%%%%%%%%%%%%%%

%%%%%%%%%%%%%%%%%%%%%
\begin{example}
\label{ex3:hjts}
Let $A$ be an associative algebra and $V = M_{pq}(A)$ be the $\bk$-module of $p \times q$ matrices with entries in $A$.  Then $(V,\{,,\})$ is a Jordan triple system with triple product
\[
\{xyz\} = xy^tz + zy^tx
\]
for all $x,y,z \in V$, where $y^t$ is the transpose of the matrix $y$ and the products on the right-hand side denote matrix multiplication.  Let $\alpha \colon A \to A$ be an algebra morphism.  Then $\alpha$ extends entry-wise to a linear map $\alpha \colon V \to V$ that is compatible with matrix multiplication (whenever it is defined) and transpose.  So $\alpha$ is a morphism of Jordan triple systems.  By Corollary \ref{cor2:twisttriple} we have a multiplicative Hom-Jordan triple system $V_\alpha = (V,\{,,\}_\alpha,\alpha)$, where $\{,,\}_\alpha = \alpha\circ\{,,\}$.
\qed
\end{example}
%%%%%%%%%%%%%%%%%%%%%

%%%%%%%%%%%%%%%%%%%%%
\begin{example}
\label{ex4:hjts}
Let $A$ be an associative algebra and $\theta$ be an involutive anti-morphism on $A$, i.e., $\theta$ is a linear automorphism on $A$ satisfying $\theta^2 = Id$ and $\theta(ab) = \theta(b)\theta(a)$ for all $a, b \in A$.  Then $(A,\{,,\})$ is a Jordan triple system with triple product
\[
\{xyz\} = x\theta(y)z + z \theta(y)x.
\]
If $\alpha \colon A \to A$ is an algebra morphism commuting with $\theta$, then $\alpha$ is a morphism of Jordan triple systems on $(A,\{,,\})$.  By Corollary \ref{cor2:twisttriple} we have a multiplicative Hom-Jordan triple system $A_\alpha = (A,\{,,\}_\alpha,\alpha)$, where $\{,,\}_\alpha = \alpha\circ\{,,\}$.
\qed
\end{example}
%%%%%%%%%%%%%%%%%%%%%

Recall that a \textbf{Jordan algebra} \cite{albert,jacobson2,schafer} $(A,\mu)$ has a bilinear multiplication $\mu$ that is commutative and that satisfies the \emph{Jordan identity}
\begin{equation}
\label{jordanid}
(x^2y)x = x^2(yx)
\end{equation}
for all $x,y \in A$, where $x^2 = \mu(x,x)$.  It is known that every Jordan algebra $(A,\mu)$ gives rise to a Jordan triple system $(A,\{,,\})$ with triple product
\begin{equation}
\label{jt}
\{xyz\} = x(yz) + (xy)z - y(xz).
\end{equation}
It can be proved using some facts about Jordan algebras that can be found in \cite{schafer} (Ch. IV).  Any algebra self-morphism on the Jordan algebra $A$ gives a morphism $(A,\{,,\}) \to (A,\{,,\})$ of Jordan triple systems.  In view of Corollary \ref{cor2:twisttriple}, we have the following method for constructing Hom-Jordan triple systems from Jordan algebras.

\begin{corollary}
\label{cor:jhjts}
Let $(A,\mu)$ be a Jordan algebra and $\alpha \colon A \to A$ be an algebra morphism.  Then $A_\alpha = (A,\{,,\}_\alpha = \alpha \circ \{,,\},\alpha)$ is a multiplicative Hom-Jordan triple system, where $\{,,\}$ is the triple product in \eqref{jt}.
\end{corollary}

The following example illustrates how Corollary \ref{cor:jhjts} can be used to construct Hom-Jordan triple systems.  Recall that an \textbf{alternative algebra} $(A,\mu)$ \cite{schafer} has a multiplication $\mu \colon A^{\otimes 2} \to A$ such that the associator
\begin{equation}
\label{associator}
as(x,y,z) = (xy)z - x(yz)
\end{equation}
is anti-symmetric.  In particular, every associative algebra is alternative.

%%%%%%%%%%%%%%%%%%%%%
\begin{example}
\label{ex1:jhjts}
The octonions $\oct$ \cite{baez,okubo,schafer} is the eight-dimensional alternative (but not associative) algebra with basis $\{e_0,\ldots,e_7\}$ and the following multiplication table:
\begin{center}
\begin{tabular}{c|c|c|c|c|c|c|c|c}
$\mu$ & $e_0$ & $e_1$ & $e_2$ & $e_3$ & $e_4$ & $e_5$ & $e_6$ & $e_7$ \\\hline
$e_0$ & $e_0$ & $e_1$ & $e_2$ & $e_3$ & $e_4$ & $e_5$ & $e_6$ & $e_7$ \\\hline
$e_1$ & $e_1$ & $-e_0$ & $e_4$ & $e_7$ & $-e_2$ & $e_6$ & $-e_5$ & $-e_3$ \\\hline
$e_2$ & $e_2$ & $-e_4$ & $-e_0$ & $e_5$ & $e_1$ & $-e_3$ & $e_7$ & $-e_6$ \\\hline
$e_3$ & $e_3$ & $-e_7$ & $-e_5$ & $-e_0$ & $e_6$ & $e_2$ & $-e_4$ & $e_1$ \\\hline
$e_4$ & $e_4$ & $e_2$ & $-e_1$ & $-e_6$ & $-e_0$ & $e_7$ & $e_3$ & $-e_5$ \\\hline
$e_5$ & $e_5$ & $-e_6$ & $e_3$ & $-e_2$ & $-e_7$ & $-e_0$ & $e_1$ & $e_4$ \\\hline
$e_6$ & $e_6$ & $e_5$ & $-e_7$ & $e_4$ & $-e_3$ & $-e_1$ & $-e_0$ & $e_2$ \\\hline
$e_7$ & $e_7$ & $e_3$ & $e_6$ & $-e_1$ & $e_5$ & $-e_4$ & $-e_2$ & $-e_0$ \\
\end{tabular}
\end{center}
For an octonion $x = \sum_{i=0}^7 b_ie_i$ with each $b_i \in \bk$, its \emph{conjugate} is defined as the octonion $\xbar = b_0e_0 - \sum_{i=1}^7 b_ie_i$.

Consider the $27$-dimensional $\bk$-module $M^8_3$ consisting of $3 \times 3$ Hermitian octonionic matrices, i.e., matrices of the form
\[
X =
\begin{pmatrix}
a_1 & x & y\\
\xbar & a_2 & z\\
\ybar & \zbar & a_3
\end{pmatrix}
\]
with each $a_i \in \bk$ and $x,y,z \in \oct$, where $a_i = a_ie_0$ for the diagonal elements.  This $\bk$-module $M^8_3$ is an exceptional simple Jordan algebra when equipped with the Jordan product
\[
X \ast Y = \frac{1}{2}(XY + YX),
\]
where $XY$ and $YX$ are the usual matrix multiplication \cite{jvw,schafer}.

Let $\alpha \colon \oct \to \oct$ be any unit-preserving and conjugate-preserving algebra morphism, i.e., $\alpha(e_0) = e_0$ and $\alpha(\xbar) = \overline{\alpha(x)}$ for all $x \in \oct$.  Then it extends entry-wise to a linear map $\alpha \colon M^8_3 \to M^8_3$, which respects matrix multiplication and hence also the Jordan product $\ast$.  By Corollary \ref{cor:jhjts} we have a multiplicative Hom-Jordan triple system
\[
(M^8_3)_\alpha = (M^8_3,\{,,\}_\alpha = \alpha \circ \{,,\},\alpha),
\]
where
\[
\{XYZ\} = X * (Y * Z) + (X * Y)* Z - Y * (X * Z)
\]
for $X,Y,Z \in M^8_3$.

There are many unit-preserving and conjugate-preserving algebra morphisms on $\oct$.  For example, the algebra automorphism $\alpha \colon \oct \to \oct$ given by
\begin{equation}
\label{octaut}
\begin{split}
\alpha(e_0) = e_0,\quad \alpha(e_1) = e_5,\quad \alpha(e_2) = e_6,\quad \alpha(e_3) = e_7,\\
\alpha(e_4) = e_1,\quad \alpha(e_5) = e_2,\quad \alpha(e_6) = e_3,\quad \alpha(e_7) = e_4
\end{split}
\end{equation}
is easily checked to be unit-preserving and conjugate-preserving.
There is a more conceptual description of this algebra automorphism on $\oct$.  Note that $e_1$ and $e_2$ anti-commute, and $e_3$ anti-commutes with $e_1$, $e_2$, and $e_1e_2 = e_4$.  Such a triple $(e_1,e_2,e_3)$ is called a \textbf{basic triple} in \cite{baez}.  Another basic triple is $(e_5,e_6,e_7)$.  Then $\alpha$ in \eqref{octaut} is the unique automorphism on $\oct$ that sends the basic triple $(e_1,e_2,e_3)$ to the basic triple $(e_5,e_6,e_7)$.
\qed
\end{example}
%%%%%%%%%%%%%%%%%%%%%

It is well-known that alternative algebras are Jordan admissible \cite{schafer}.  In other words, if $(A,\mu)$ is an alternative algebra, then the anti-commutator algebra
\[
A^+ = (A,* = (\mu + \muop)/2)
\]
is a Jordan algebra, where
\[
x * y = \frac{xy + yx}{2}.
\]
The Hom-version of this fact, that Hom-alternative algebras are Hom-Jordan admissible, is proved in \cite{yau12}, but we do not need that result here.  Every algebra self-morphism on an alternative algebra $A$ is also multiplicative with respect to the Jordan product $*$.  Therefore, we have the following consequence of Corollary \ref{cor:jhjts}, which gives a way of constructing Hom-Jordan triple systems from alternative algebras.

\begin{corollary}
\label{cor:althjts}
Let $(A,\mu)$ be an alternative algebra and $\alpha \colon A \to A$ be an algebra morphism.  Then
\[
A^+_\alpha = (A,\{,,\}_\alpha,\alpha)
\]
is a multiplicative Hom-Jordan triple system, where
\[
\{xyz\}_\alpha = \alpha\left(x * (y * z) + (x * y) * z - y * (x * z)\right)
\]
and $x*y = (xy + yx)/2$.
\end{corollary}

We close this section with some examples of Hom-Lie triple systems, which are automatically ternary Hom-Nambu algebras.  Note that the ternary Nambu algebras $(V,[,,])$ in Examples \ref{ex:okubo} and \ref{ex:fermion} are actually Lie triple systems, so in each case $V_\alpha = (V,[,,]_\alpha=\alpha\circ[,,],\alpha)$ is a multiplicative Hom-Lie triple system.

% Maltsev to LTS (Loos) to Hom-LTS
In \cite{loos} Loos showed that every Maltsev algebra has an underlying Lie triple system.  A \textbf{Maltsev algebra} \cite{maltsev,sagle} $(A,\mu)$ has an anti-symmetric multiplication $\mu \colon A^{\otimes 2} \to A$ that satisfies the \emph{Maltsev identity}
\begin{equation}
\label{maltsev}
J'(x,y,xz) = J'(x,y,z)x
\end{equation}
for all $x,y,z \in A$, where $\mu(x,y) = xy$ and
\[
J'(x,y,z) = (xy)z + (zx)y + (yz)x
\]
is the Jacobian.  In particular, Lie algebras are examples of Maltsev algebras.   Maltsev algebras play an important role in the geometry of smooth loops.   Just as the tangent algebra of a Lie group is a Lie algebra, the tangent algebra of a locally analytic Moufang loop is a Maltsev algebra \cite{kerdman,kuzmin,maltsev,nagy,sabinin}.  The Hom-version of Maltsev(-admissible) algebras are studied in \cite{yau12}.

According to a result in \cite{loos}, if $(A,\mu)$ is a Maltsev algebra, then $(A,[,,])$ is a Lie triple system with triple product
\begin{equation}
\label{loosproduct}
[xyz] = 2(xy)z - (zx)y - (yz)x.
\end{equation}
Any algebra self-morphism on $A$ is also a morphism of triple systems.  In view of Corollary \ref{cor2:twisttriple}, we have the following method of constructing Hom-Lie triple systems from Maltsev algebras.

\begin{corollary}
\label{cor:maltsevhlts}
Let $(A,\mu)$ be a Maltsev algebra and $\alpha \colon A \to A$ be an algebra morphism.  Then $A_\alpha = (A,[,,]_\alpha = \alpha \circ [,,],\alpha)$ is a multiplicative Hom-Lie triple system, where $[,,]$ is the triple product in \eqref{loosproduct}.
\end{corollary}

% alternative to Maltsev (Bruck-Kleinfeld) to LTS to Hom-LTS
Maltsev observed in \cite{maltsev} (see also \cite{bk}) that every alternative algebra $(A,\mu)$ is Maltsev-admissible, i.e., the commutator algebra $A^- = (A,[,] = \mu - \muop)$ is a Maltsev algebra.  The Hom-version of this fact, that Hom-alternative algebras are Hom-Maltsev admissible, is proved in \cite{yau12}, but we do not need that result here.  Every algebra morphism on $(A,\mu)$ is also an algebra morphism on $A^-$.  Using Corollary \ref{cor:maltsevhlts}, we thus have the following method of constructing Hom-Lie triple systems from alternative algebras.

\begin{corollary}
\label{cor:althlts}
Let $(A,\mu)$ be an alternative algebra and $\alpha \colon A \to A$ be an algebra morphism.  Then $A^-_\alpha = (A,[,,]_\alpha,\alpha)$ is a multiplicative Hom-Lie triple system, where
\[
[xyz]_\alpha = \alpha\left(2[[x,y],z] - [[z,x],y] - [[y,z],x]\right)
\]
and $[,] = \mu - \muop$ is the commutator bracket of $\mu$.
\end{corollary}

%%%%%%%%%%%%%%%%%%%%%%%%%%%%%%%%%%%%%%%%%%%%
\section{Hom-Lie triple systems}
\label{sec:meyberg}
%%%%%%%%%%%%%%%%%%%%%%%%%%%%%%%%%%%%%%%%%%%%

The purpose of this section is to show that Hom-Lie triple systems (and hence ternary Hom-Nambu algebras) arise from Hom-Jordan triple systems, ternary totally Hom-associative algebras, Hom-associative algebras, and Hom-Lie algebras.

Here is the first main result of this section, which says that every ternary totally Hom-associative algebra with equal twisting maps has an underlying Hom-Jordan triple system.

% hom-ats to hom-jts
\begin{theorem}
\label{thm:atsjts}
Let $(A,(,,),\alpha)$ be a ternary totally Hom-associative algebra whose twisting maps are equal.  Define the triple product
\begin{equation}
\label{hatshjts}
\{xyz\} = (xyz) + (zyx)
\end{equation}
for $x,y,z \in A$.  Then $J(A) = (A,\{,,\},\alpha)$ is a Hom-Jordan triple system.  Moreover, if $A$ is multiplicative, then so is $J(A)$.
\end{theorem}

\begin{proof}
It is immediate from the definition that the triple product $\{,,\}$ in \eqref{hatshjts} is symmetric in $x$ and $z$.  To check the Hom-Jordan triple identity \eqref{homjtsid}, note that the left-hand side of \eqref{homjtsid} is:
\begin{equation}
\label{hjtsl}
\begin{split}
&\{\alpha(x)\alpha(y)\{uvw\}\} - \{\alpha(u)\alpha(v)\{xyw\}\}\\
&= (\alpha(x)\alpha(y)(uvw)) + (\alpha(x)\alpha(y)(wvu)) + ((uvw)\alpha(y)\alpha(x)) + ((wvu)\alpha(y)\alpha(x))\\
&\relphantom{} - (\alpha(u)\alpha(v)(xyw)) - (\alpha(u)\alpha(v)(wyx)) - ((xyw)\alpha(v)\alpha(u)) - ((wyx)\alpha(v)\alpha(u))\\
&=  (\alpha(x)\alpha(y)(uvw)) + ((wvu)\alpha(y)\alpha(x)) - (\alpha(u)\alpha(v)(xyw)) - ((wyx)\alpha(v)\alpha(u)).
\end{split}
\end{equation}
In the last equality above, we used ternary Hom-associativity \eqref{thomass} twice.  Likewise, the right-hand side of \eqref{homjtsid} is:
\begin{equation}
\label{hjtsr}
\begin{split}
&\{\{xyu\}\alpha(v)\alpha(w)\} - \{\alpha(u)\{yxv\}\alpha(w)\}\\
&= ((xyu)\alpha(v)\alpha(w)) + ((uyx)\alpha(v)\alpha(w)) + (\alpha(w)\alpha(v)(xyu)) + (\alpha(w)\alpha(v)(uyx))\\
&\relphantom{} - (\alpha(u)(yxv)\alpha(w)) - (\alpha(u)(vxy)\alpha(w)) - (\alpha(w)(yxv)\alpha(u)) - (\alpha(w)(vxy)\alpha(u))\\
&= ((xyu)\alpha(v)\alpha(w)) + (\alpha(w)\alpha(v)(uyx)) - (\alpha(u)(vxy)\alpha(w)) - (\alpha(w)(yxv)\alpha(u)).
\end{split}
\end{equation}
Using ternary Hom-associativity \eqref{thomass} four times, one observes that \eqref{hjtsl} and \eqref{hjtsr} are equal, showing that $J(A)$ is a Hom-Jordan triple system.  The second assertion regarding multiplicativity is immediate from the definition.
\end{proof}

The next result says that every Hom-Jordan triple system with equal twisting maps has an underlying Hom-Lie triple system, and hence a ternary Hom-Nambu algebra.  It generalizes an observation of Meyberg \cite{meyberg2} (XI Theorem I) that a Jordan triple system gives rise to a Lie triple system.

% hom-jts to hom-lts
\begin{theorem}
\label{thm:jtslts}
Let $(J,\{,,\},\alpha)$ be a Hom-Jordan triple system with equal twisting maps.  Define the triple product
\begin{equation}
\label{hjtshlts}
[xyz] = \{xyz\} - \{yxz\}
\end{equation}
for $x,y,z \in J$.  Then $L(J) = (J,[,,],\alpha)$ is a Hom-Lie triple system.  Moreover, if $J$ is multiplicative, then so is $L(J)$.
\end{theorem}

\begin{proof}
Both the left anti-symmetry and the ternary Jacobi identity for $[,,]$ are immediate from the definition \eqref{hjtshlts}.  To check the ternary Hom-Nambu identity \eqref{homnambu}, note that the left-hand side of \eqref{homnambu} is:
\begin{equation}
\label{homnambul}
\begin{split}
&[\alpha(x)\alpha(y)[uvw]]\\
&= \{\alpha(x)\alpha(y)\{uvw\}\} - \{\alpha(x)\alpha(y)\{vuw\}\} - \{\alpha(y)\alpha(x)\{uvw\}\} + \{\alpha(y)\alpha(x)\{vuw\}\}\\
&= \{\alpha(u)\alpha(v)\{xyw\}\} + \{\{xyu\}\alpha(v)\alpha(w)\} - \{\alpha(u)\{yxv\}\alpha(w)\}\\
&\relphantom{} - \{\alpha(v)\alpha(u)\{xyw\}\} - \{\{xyv\}\alpha(u)\alpha(w)\} + \{\alpha(v)\{yxu\}\alpha(w)\}\\
&\relphantom{} - \{\alpha(u)\alpha(v)\{yxw\}\} - \{\{yxu\}\alpha(v)\alpha(w)\} + \{\alpha(u)\{xyv\}\alpha(w)\}\\
&\relphantom{} + \{\alpha(v)\alpha(u)\{yxw\}\} + \{\{yxv\}\alpha(u)\alpha(w)\} - \{\alpha(v)\{xyu\}\alpha(w)\}.
\end{split}
\end{equation}
In the last equality above, we used the Hom-Jordan triple identity \eqref{homjtsid} four times.  On the other hand, expanding the right-hand side of the ternary Hom-Nambu identity \eqref{homnambu} using the definition \eqref{hjtshlts}, we obtain a sum of twelve terms.  They are exactly the twelve terms in \eqref{homnambul}, showing that $L(J)$ is a Hom-Lie triple system. The last assertion regarding multiplicativity is immediate from the definition \eqref{hjtshlts}.
\end{proof}

Theorem \ref{thm:jtslts} can be applied to the various Hom-Jordan triple systems in section \ref{sec:homtriple} to obtain Hom-Lie triple systems.

Combining Theorems \ref{thm:atsjts} and \ref{thm:jtslts}, we obtain the following way of constructing a Hom-Lie triple system, and hence a ternary Hom-Nambu algebra, from a ternary totally Hom-associative algebra.

% hom-ats to hom-lts
\begin{corollary}
\label{cor:atslts}
Let $(A,(,,),\alpha)$ be a ternary totally Hom-associative algebra with equal twisting maps.  Define the triple product
\begin{equation}
\label{hatshlts}
[xyz] = (xyz) - (yxz) - (zxy) + (zyx)
\end{equation}
for $x,y,z \in A$.  Then $L(A) = (A,[,,],\alpha)$ is a Hom-Lie triple system.  Moreover, if $A$ is multiplicative, then so is $L(A)$.
\end{corollary}

In the special case $\alpha = Id$, Corollary \ref{cor:atslts} becomes the following result.

\begin{corollary}
\label{cor2:atslts}
Let $(A,(,,))$ be a ternary ring.  Then $(A,[,,])$ is a Lie triple system, where $[,,]$ is the triple product in \eqref{hatshlts}.
\end{corollary}

Next we discuss how ternary totally Hom-associative algebras and Hom-Lie triple systems arise from Hom-associative algebras and Hom-Lie algebras.

\begin{definition}
\label{def:homalg}
\begin{enumerate}
\item
By a \textbf{Hom-algebra} we mean a binary Hom-algebra $(A,\mu,\alpha)$.  Multiplicativity, (weak) morphisms, and anti-symmetry are defined as in Definition \ref{def:nhomalgebra}.
\item
The \textbf{Hom-associator} of a Hom-algebra $(A,\mu,\alpha)$ is the trilinear map $as_A \colon A^{\otimes 3} \to A$ defined as
\begin{equation}
\label{homassociator}
as_A = \mu \circ (\mu \otimes \alpha - \alpha \otimes \mu).
\end{equation}
A \textbf{Hom-associative algebra} \cite{ms} is a Hom-algebra whose Hom-associator is equal to $0$.
\item
A \textbf{Hom-Lie algebra} \cite{hls,ms} is a binary Hom-Nambu-Lie algebra $(L,[,],\alpha)$.  In this case, the binary Hom-Nambu identity $J^2_L = 0$ is called the \textbf{Hom-Jacobi identity}, which by anti-symmetry is equivalent to
\begin{equation}
\label{homjacobi}
[[x,y],\alpha(z)] + [[z,x],\alpha(y)] + [[y,z],\alpha(x)] = 0
\end{equation}
for all $x,y,z \in L$.
\end{enumerate}
\end{definition}

Construction results along the lines of Corollary \ref{cor3:twist} and examples of Hom-associative and Hom-Lie algebras can be found in \cite{yau2}.

The following result says that every Hom-associative algebra has an underlying ternary totally Hom-associative algebra.

% hom-associative algebras to hom-ats
\begin{theorem}
\label{thm:haats}
Let $(A,\mu,\alpha)$ be a Hom-associative algebra.  Then
\[
A_T = (A,(,,) = \mu \circ (\mu \otimes \alpha),\alpha^2)
\]
is a ternary totally Hom-associative algebra.  Moreover, if $A$ is multiplicative, then so is $A_T$.
\end{theorem}

\begin{proof}
If $\alpha$ is multiplicative with respect to $\mu$, then it is clear that $\alpha^2$ is multiplicative with respect to the triple product $(,,)$.  To check ternary Hom-associativity \eqref{thomass} for $A_T$, note that $as_A = 0$ means that $\mu \circ (\mu \otimes \alpha) = \mu \circ (\alpha \otimes \mu)$.  Using this repeatedly, we compute as follows, where $\mu(x,y)$ is written as the juxtaposition $xy$:
\[
\begin{split}
((uvw)\alpha^2(x)\alpha^2(y))
&= (((uv)\alpha(w))\alpha^2(x))\alpha^3(y)\\
&= ((\alpha(u)(vw))\alpha^2(x))\alpha^3(y)\\
&= (\alpha^2(u)((vw)\alpha(x)))\alpha^3(y)\\
&= \alpha^3(u)(((vw)\alpha(x))\alpha^2(y))\\
&= (\alpha^2(u)(vwx)\alpha^2(y))\\
&= \alpha^3(u)((\alpha(v)(wx))\alpha^2(y))\\
&= \alpha^3(u)((\alpha^2(v)((wx)\alpha(y))))\\
&= (\alpha^2(u)\alpha^2(v)(wxy)).
\end{split}
\]
This shows that $A_T$ is a ternary totally Hom-associative algebra.
\end{proof}

The following result is the Hom-Lie version of Theorem \ref{thm:haats}.  It says that every multiplicative Hom-Lie algebra has an underlying multiplicative Hom-Lie triple system, and hence a multiplicative ternary Hom-Nambu algebra.

% hom-lie algebras to hom-lts
\begin{theorem}
\label{thm:hllts}
Let $(L,[,],\alpha)$ be a multiplicative Hom-Lie algebra.  Then
\[
L_T = (L,[,,] = [,] \circ ([,] \otimes \alpha),\alpha^2)
\]
is a multiplicative Hom-Lie triple system.
\end{theorem}

\begin{proof}
The left-antisymmetry of $[,,]$ follows from the anti-symmetry of $[,]$. The ternary Jacobi identity \eqref{homltsid} for $[,,]$ is exactly the Hom-Jacobi identity \eqref{homjacobi} for $[,]$.  The ternary Hom-Nambu identity \eqref{homnambu} for $[,,]$ is proved by using the multiplicativity of $\alpha$ and the Hom-Jacobi identity repeatedly.  Indeed, we have:
\[
\begin{split}
&[\alpha^2(x)\alpha^2(y)[uvw]]\\
&= [\alpha[\alpha(x),\alpha(y)],[[\alpha(u),\alpha(v)],\alpha^2(w)]]\\
&= [[[\alpha(x),\alpha(y)],[\alpha(u),\alpha(v)]],\alpha^3(w)] + [\alpha[\alpha(u),\alpha(v)],[[\alpha(x),\alpha(y)],\alpha^2(w)]]\\
&= [[[[x,y],\alpha(u)],\alpha^2(v)],\alpha^3(w)] + [[\alpha^2(u),[[x,y],\alpha(v)]],\alpha^3(w)]\\
&\relphantom{} + [[\alpha^2(u),\alpha^2(v)],\alpha[[x,y],\alpha(w)]]\\
&= [[xyu]\alpha^2(v)\alpha^2(w)] + [\alpha^2(u)[xyv]\alpha^2(w)] + [\alpha^2(u)\alpha^2(v)[xyw]].
\end{split}
\]
This shows that $L_T$ satisfies the ternary Hom-Nambu identity \eqref{homnambu}.
\end{proof}

Note that in Theorem \ref{thm:haats}, a Hom-associative algebra gives rise to a ternary totally Hom-associative algebra, even if it is not multiplicative.  In contrast, in Theorem \ref{thm:hllts} the multiplicativity assumption is necessary.

Combining Corollary \ref{cor:atslts} and Theorem \ref{thm:haats}, we obtain the following method of constructing Hom-Lie triple systems, and hence ternary Hom-Nambu algebras, from Hom-associative algebras.

\begin{corollary}
\label{cor:hahnambu}
Let $(A,\mu,\alpha)$ be a Hom-associative algebra.  Then
\[
A_L = (A,[,,],\alpha^2)
\]
is a Hom-Lie triple system, where (writing $\mu(x,y)$ as $xy$)
\[
[xyz] = (xy)\alpha(z) - (yx)\alpha(z) - (zx)\alpha(y) + (zy)\alpha(x)
\]
for $x,y,z\in A$.  Moreover, if $A$ is multiplicative, then so is $A_L$.
\end{corollary}

Observe that the Hom-Lie triple systems in Corollary \ref{cor:atslts}, Theorem \ref{thm:hllts}, and Corollary \ref{cor:hahnambu} are usually not anti-symmetric.  In particular, they are ternary Hom-Nambu algebras that are usually not ternary Hom-Nambu-Lie algebras.

%%%%%%%%%%%%%%%%%%%%%%%%%%%%%%%%%%%%%%%%%%%%%%
\section{Hom-Nambu algebras of higher arities}
\label{sec:arity}
%%%%%%%%%%%%%%%%%%%%%%%%%%%%%%%%%%%%%%%%%%%%%%

The purpose of this section is to observe that every multiplicative $n$-ary Hom-Nambu algebra (Definition \ref{def:homnambulie}) gives rise to a sequence of multiplicative Hom-Nambu algebras of increasingly higher arities.  This result is a consequence of the following result.

\begin{theorem}
\label{thm:higher}
Let $(L,\bracket,\alpha)$ be a multiplicative $n$-ary Hom-Nambu algebra. Define the $(2n-1)$-ary product
\[
[x_1,\ldots,x_{2n-1}]^{(1)} = [[x_1,\ldots,x_n],\alpha(x_{n+1}),\ldots,\alpha(x_{2n-1})]
\]
for $x_i \in L$.  Then
\[
L^1 = (L,\bracket^{(1)},\alpha^2)
\]
is a multiplicative $(2n-1)$-ary Hom-Nambu algebra.
\end{theorem}

\begin{proof}
The multiplicativity of $L^1$ follows from that of $L$.  The proof of the $(2n-1)$-ary Hom-Nambu identity for $L^1$, $J^{2n-1}_{L^1} = 0$, is a slight generalization of the proof of Theorem \ref{thm:hllts}, using the multiplicativity of $\alpha$ and the $n$-ary Hom-Nambu identity $J^n_L = 0$ twice.  To prove it, we use the abbreviations in \eqref{xij} and the shorthand
\[
\begin{split}
[x_1,\ldots, x_{2n-2},y]^{(1)} &= [x,y]^{(1)}\\
&= [[x_{1,n}],\alpha(x_{n+1,2n-2}),\alpha(y)]
\end{split}
\]
for $x_1,\ldots,x_{2n-2},y \in L$.  In particular, we have
\[
\alpha[x,y]^{(1)} = [\alpha[x_{1,n}],\alpha^2(x_{n+1,2n-2}),\alpha^2(y)]
\]
by multiplicativity.

The $(2n-1)$-ary Hom-Jacobian \eqref{homjacobian} of $L^1$ is
\begin{equation}
\label{j2n1}
\begin{split}
&[\alpha^2(x_{1,2n-2}),[y_{1,2n-1}]^{(1)}]^{(1)}
- \sum_{i=1}^{2n-1} [\alpha^2(y_{1,i-1}),[x,y_i]^{(1)},\alpha^2(y_{i+1,2n-1})]^{(1)}\\
&= [[\alpha^2(x_{1,n})],\alpha^3(x_{n+1,2n-2}), \alpha[[y_{1,n}],\alpha(y_{n+1,2n-1})]]\\
&\relphantom{} - \sum_{i=1}^n [[\alpha^2(y_{1,i-1}),[x,y_i]^{(1)},\alpha^2(y_{i+1,n})], \alpha^3(y_{n+1,2n-1})]\\
&\relphantom{} - \sum_{i=n+1}^{2n-1} [[\alpha^2(y_{1,n})],\alpha^3(y_{n+1,i-1}), \alpha[x,y_i]^{(1)}, \alpha^3(y_{i+1,2n-1})].
\end{split}
\end{equation}
Using the $n$-ary Hom-Nambu identity $J^n_L = 0$ and multiplicativity repeatedly, the first term on the right-hand side of \eqref{j2n1} becomes
\[
\begin{split}
& [\alpha^2[x_{1,n}],\alpha^3(x_{n+1,2n-2}), [[\alpha(y_{1,n})],\alpha^2(y_{n+1,2n-1})]]\\
&= [[\alpha[x_{1,n}],\alpha^2(x_{n+1,2n-2}),[\alpha(y_{1,n})]], \alpha^3(y_{n+1,2n-1})]\\
&\relphantom{} + \sum_{j=n+1}^{2n-1} [\alpha[\alpha(y_{1,n})],\alpha^3(y_{n+1,j-1}), [\alpha[x_{1,n}], \alpha^2(x_{n+1,2n-2}), \alpha^2(y_j)], \alpha^3(y_{j+1,2n-1})]\\
&= \sum_{j=1}^n [[\alpha^2(y_{1,j-1}),[x,y_j]^{(1)},\alpha^2(y_{j+1,n})], \alpha^3(y_{n+1,2n-1})]\\
&\relphantom{} + \sum_{j=n+1}^{2n-1} [[\alpha^2(y_{1,n})],\alpha^3(y_{n+1,j-1}), \alpha[x,y_j]^{(1)}, \alpha^3(y_{j+1,2n-1})].
\end{split}
\]
These two sums cancel out with the last two sums in \eqref{j2n1}.  This shows that the $(2n-1)$-ary Hom-Jacobian \eqref{homjacobian} of $L^1$ is equal to $0$, as desired.
\end{proof}

It is easy to see that the construction in Theorem \ref{thm:higher} gives a functor from the category of multiplicative $n$-ary Hom-Nambu algebras (where the maps are the morphisms) to the category of multiplicative $(2n-1)$-ary Hom-Nambu algebras.   Iterating Theorem \ref{thm:higher}, we obtain a sequence of Hom-Nambu algebras of exponentially higher arities, as in the following result.  We use the abbreviations in \eqref{xij}.

\begin{corollary}
\label{cor1:higher}
Let $(L,\bracket,\alpha)$ be a multiplicative $n$-ary Hom-Nambu algebra. For $k \geq 0$ define the $(2^k(n-1)+1)$-ary product $\bracket^{(k)}$ inductively by setting $\bracket^{(0)} = \bracket$ and
\[
[x_{1,2^k(n-1)+1}]^{(k)} = [[x_{1,2^{k-1}(n-1)+1}]^{(k-1)}, \alpha^{2^{k-1}}(x_{2^{k-1}(n-1)+2,2^k(n-1)+1})]^{(k-1)}
\]
for $k \geq 1$ and $x_i \in L$.  Then
\[
L^k = (L,\bracket^{(k)},\alpha^{2^k})
\]
is a multiplicative $(2^k(n-1)+1)$-ary Hom-Nambu algebra for each $k \geq 0$.
\end{corollary}

The following result is the special case of Corollary \ref{cor1:higher} when $\alpha$ is the identity map.

\begin{corollary}
\label{cor2:higher}
Let $(L,\bracket)$ be an $n$-ary Nambu algebra.  For $k \geq 0$ define the $(2^k(n-1)+1)$-ary product $\bracket^{(k)}$ inductively by setting $\bracket^{(0)} = \bracket$ and
\[
[x_{1,2^k(n-1)+1}]^{(k)} = [[x_{1,2^{k-1}(n-1)+1}]^{(k-1)}, x_{2^{k-1}(n-1)+2,2^k(n-1)+1}]^{(k-1)}
\]
for $k \geq 1$ and $x_i \in L$.  Then $L^k = (L,\bracket^{(k)})$ is a $(2^k(n-1)+1)$-ary Nambu algebra for each $k \geq 0$.
\end{corollary}

We can apply Corollary \ref{cor1:higher} to any of the multiplicative $n$-ary Hom-Nambu algebras from the previous sections, such as the ones in Corollary \ref{cor2:twist} and Corollary \ref{cor1.5:twisttriple}.  (Recall that a Hom-Lie triple system (Definition \ref{def:hjts}) is automatically a ternary Hom-Nambu algebra.)

If $L$ is a multiplicative binary Hom-Nambu algebra, such as a multiplicative Hom-Lie algebra, then $L^k$ in Corollary \ref{cor1:higher} is a multiplicative $(2^k+1)$-ary Hom-Nambu algebra.

If $L$ is a multiplicative ternary Hom-Nambu algebra, then $L^k$ in Corollary \ref{cor1:higher} is a multiplicative $(2^{k+1}+1)$-ary Hom-Nambu algebra.  For example, using Corollary \ref{cor1:higher} with Theorem \ref{thm:jtslts}, we obtain a sequence of multiplicative $(2^{k+1}+1)$-ary Hom-Nambu algebras from any multiplicative Hom-Jordan triple system.  Let us discuss the cases $k=1$ and $2$ in the following example.

%%%%%%%%%%%%%%%%%%%%%
\begin{example}
\label{ex:hjtsnambu}
Let $(J,\{,,\},\alpha)$ be a multiplicative Hom-Jordan triple system and $L(J) = (J,[,,],\alpha)$ be the multiplicative Hom-Lie triple system (hence multiplicative ternary Hom-Nambu algebra) in Theorem \ref{thm:jtslts}, where
\[
[xyz] = \{xyz\} - \{yxz\}.
\]
By Corollary \ref{cor1:higher} there is a multiplicative $5$-ary Hom-Nambu algebra
\[
L(J)^1 = (J,\bracket^{(1)},\alpha^2)
\]
with
\[
\begin{split}
[x_1,\ldots,x_5]^{(1)}
&= [[x_1,x_2,x_3],\alpha(x_4),\alpha(x_5)]\\
&= \{\{x_1,x_2,x_3\},\alpha(x_4),\alpha(x_5)\} - \{\alpha(x_4),\{x_1,x_2,x_3\},\alpha(x_5)\}\\
&\relphantom{} - \{\{x_2,x_1,x_3\},\alpha(x_4),\alpha(x_5)\} + \{\alpha(x_4),\{x_2,x_1,x_3\},\alpha(x_5)\}.
\end{split}
\]
By Corollary \ref{cor1:higher} again, there is a multiplicative $9$-ary Hom-Nambu algebra
\[
L(J)^2 = (J,\bracket^{(2)},\alpha^4)
\]
with
\begin{equation}
\label{2bracket}
\begin{split}
[x_1,\ldots,x_9]^{(2)}
&= [[x_1,\ldots,x_5]^{(1)},\alpha^2(x_6),\ldots,\alpha^2(x_9)]^{(1)}\\
&= [[[x_{1,3}],\alpha(x_{4,5})],\alpha^2(x_{6,9})]^{(1)}\\
&= [[[[x_{1,3}],\alpha(x_{4,5})],\alpha^2(x_{6,7})],\alpha^3(x_{8,9})].
\end{split}
\end{equation}
One can write this last expression in terms of the Hom-Jordan triple product $\{,,\}$.  It involves $16$ terms, each one involving a composition of four copies of $\{,,\}$.
\qed
\end{example}
%%%%%%%%%%%%%%%%%%%%%

On the other hand, using Corollary \ref{cor1:higher} with Corollary \ref{cor:atslts}, we obtain a sequence of multiplicative $(2^{k+1}+1)$-ary Hom-Nambu algebras from any multiplicative ternary totally Hom-associative algebra.  For example, if $(A,\mu,\alpha)$ is a multiplicative Hom-associative algebra and $A_T$ is the multiplicative ternary totally Hom-associative algebra in Theorem \ref{thm:haats}, then Corollary \ref{cor1:higher} and Corollary \ref{cor:atslts} yield a sequence of multiplicative $(2^{k+1}+1)$-ary Hom-Nambu algebras.  Let us discuss the cases $k=1$ and $2$ in the following example.

%%%%%%%%%%%%%%%%%%%%
\begin{example}
\label{ex:hatsnambu}
Let $(A,(,,),\alpha)$ be a multiplicative ternary totally Hom-associative algebra and $L(A) = (A,[,,],\alpha)$ be the multiplicative Hom-Lie triple system (hence multiplicative ternary Hom-Nambu algebra) in Corollary \ref{cor:atslts}, where
\[
[xyz] = (xyz) - (yxz) - (zxy) + (zyx).
\]
By Corollary \ref{cor1:higher} there is a multiplicative $5$-ary Hom-Nambu algebra
\[
L(A)^1 = (A,\bracket^{(1)},\alpha^2)
\]
with
\[
\begin{split}
[x_1,\ldots,x_5]^{(1)}
&= [(x_1,x_2,x_3),\alpha(x_4),\alpha(x_5)] - [(x_2,x_1,x_3),\alpha(x_4),\alpha(x_5)]\\
&\relphantom{} - [(x_3,x_1,x_2),\alpha(x_4),\alpha(x_5)] + [(x_3,x_2,x_1),\alpha(x_4),\alpha(x_5)]\\
&= ((x_1,x_2,x_3),\alpha(x_4),\alpha(x_5)) - (\alpha(x_4),(x_1,x_2,x_3),\alpha(x_5))\\
&\relphantom{} - (\alpha(x_5),(x_1,x_2,x_3),\alpha(x_4)) + (\alpha(x_5),\alpha(x_4),(x_1,x_2,x_3))\\
&\relphantom{} - ((x_2,x_1,x_3),\alpha(x_4),\alpha(x_5)) + (\alpha(x_4),(x_2,x_1,x_3),\alpha(x_5))\\
&\relphantom{} + (\alpha(x_5),(x_2,x_1,x_3),\alpha(x_4)) - (\alpha(x_5),\alpha(x_4),(x_2,x_1,x_3))\\
&\relphantom{} - ((x_3,x_1,x_2),\alpha(x_4),\alpha(x_5)) + (\alpha(x_4),(x_3,x_1,x_2),\alpha(x_5))\\
&\relphantom{} + (\alpha(x_5),(x_3,x_1,x_2),\alpha(x_4)) - (\alpha(x_5),\alpha(x_4),(x_3,x_1,x_2))\\
&\relphantom{} + ((x_3,x_2,x_1),\alpha(x_4),\alpha(x_5)) - (\alpha(x_4),(x_3,x_2,x_1),\alpha(x_5))\\
&\relphantom{} - (\alpha(x_5),(x_3,x_2,x_1),\alpha(x_4)) + (\alpha(x_5),\alpha(x_4),(x_3,x_2,x_1)).
\end{split}
\]
By Corollary \ref{cor1:higher} again, there is a multiplicative $9$-ary Hom-Nambu algebra
\[
L(A)^2 = (A,\bracket^{(2)},\alpha^4)
\]
with
\[
[x_1,\ldots,x_9]^{(2)}
= [[[[x_{1,3}],\alpha(x_{4,5})],\alpha^2(x_{6,7})],\alpha^3(x_{8,9})]
\]
as in \eqref{2bracket}.  One can write this last expression in terms of the ternary totally Hom-associative product $(,,)$.  It involves $4^4$ terms, each one involving a composition of four copies of $(,,)$.
\qed
\end{example}
%%%%%%%%%%%%%%%%%%%%

%%%%%%%%%%%%%%%%%%%%%%%%%%%%%%%%%%%%%%%%%%%%%%%%%%%
\section{Hom-Nambu(-Lie) algebras of lower arities}
\label{sec:lower}
%%%%%%%%%%%%%%%%%%%%%%%%%%%%%%%%%%%%%%%%%%%%%%%%%%%

The purpose of this section is to observe that, under suitable assumptions, an $n$-ary Hom-Nambu(-Lie) algebra with $n \geq 3$ reduces to an $(n-1)$-ary Hom-Nambu(-Lie) algebra.

First we need the following observations regarding a reduced product.  The following Lemma is inspired by \cite{poz} (Lemma 1.2).

\begin{lemma}
\label{lem:lower}
Let $(V,\bracket,\alpha=(\alpha_1,\ldots,\alpha_{n-1}))$ be an $n$-ary Hom-algebra with $n \geq 3$.  Suppose $a \in V$ satisfies $\alpha_1(a) = a$.  Define the $(n-1)$-ary product
\[
[x_1,\ldots, x_{n-1}]' = [a,x_1,\ldots,x_{n-1}]
\]
for $x_1, \ldots , x_{n-1} \in V$, and consider the $(n-1)$-ary Hom-algebra
\[
V' = (V,\bracket',\alpha'=(\alpha_2,\ldots,\alpha_{n-1})).
\]
Then the following statements hold.
\begin{enumerate}
\item
If $\bracket$ is anti-symmetric, then so is $\bracket'$.
\item
If $V$ is multiplicative, then so is $V'$.
\item
The $(n-1)$-ary Hom-Jacobian of $V'$ \eqref{homjacobian} satisfies
\begin{equation}
\label{jv'}
\begin{split}
J^{n-1}_{V'}(x_{2,n-1};y_{2,n})
&= J^n_V(a,x_{2,n-1};a,y_{2,n})\\
&\relphantom{} + [[a,x_{2,n-1},a],\alpha_1(y_2),\ldots,\alpha_{n-1}(y_n)]
\end{split}
\end{equation}
for all $x_i,y_j \in V$.
\end{enumerate}
\end{lemma}

\begin{proof}
The first assertion is clear, and it does not even require the assumption $\alpha_1(a) = a$.  The second assertion follows from the multiplicativity of $\alpha$ and the assumption $\alpha_1(a) = a$.  The equality \eqref{jv'} is obtained from the definition \eqref{homjacobian} of the $n$-ary Hom-Jacobian of $V$ by specifying $x_1 = y_1 = a$ and using the assumption $\alpha_1(a) = a$.
\end{proof}

The following result is an immediate consequence of Lemma \ref{lem:lower}.

\begin{theorem}
\label{thm:lower}
Let $(L,\bracket,\alpha=(\alpha_1,\ldots,\alpha_{n-1}))$ be an $n$-ary Hom-Nambu algebra with $n \geq 3$.  Suppose $a \in L$ satisfies
\[
\alpha_1(a) = a\quad\text{and}\quad
[a,x_{2,n-1},a] = 0\quad\text{for all $x_i \in L$}.
\]
Then
\begin{equation}
\label{l'}
L' = (L,\bracket',\alpha'=(\alpha_2,\ldots,\alpha_{n-1}))
\end{equation}
is an $(n-1)$-ary Hom-Nambu algebra, where
\[
[x_1,\ldots,x_{n-1}]' = [a,x_1,\ldots,x_{n-1}]
\]
for $x_i \in L$.  Moreover, if $L$ is multiplicative, then so is $L'$.
\end{theorem}

The assumption $[a,x_{2,n-1},a] = 0$ in Theorem \ref{thm:lower} is automatically satisfied if the $n$-ary product $\bracket$ is anti-symmetric, as in a Hom-Nambu-Lie algebra.  Therefore, we have the following consequence of Lemma \ref{lem:lower} and Theorem \ref{thm:lower}.

\begin{corollary}
\label{cor1:lower}
Let $(L,\bracket,\alpha=(\alpha_1,\ldots,\alpha_{n-1}))$ be an $n$-ary Hom-Nambu-Lie algebra with $n \geq 3$.  Suppose $a \in L$ satisfies $\alpha_1(a) = a$.  Then $L'$ in \eqref{l'} is an $(n-1)$-ary Hom-Nambu-Lie algebra.  Moreover, if $L$ is multiplicative, then so is $L'$.
\end{corollary}

If we apply Theorem \ref{thm:lower} or Corollary \ref{cor1:lower} repeatedly, then we obtain the following two results.  They tell us how to go from $n$-ary Hom-Nambu(-Lie) algebras to $(n-k)$-ary Hom-Nambu(-Lie) algebras.

\begin{corollary}
\label{cor2:lower}
Let $(L,\bracket,\alpha=(\alpha_1,\ldots,\alpha_{n-1}))$ be an $n$-ary Hom-Nambu algebra with $n \geq 3$.  Suppose for some $k \in \{1,\ldots,n-2\}$ there exist $a_i \in L$ for $1 \leq i \leq k$ satisfying
\[
\alpha_i(a_i) = a_i \quad \text{for $1 \leq i \leq k$}
\]
and
\[
[a_1,\ldots,a_j,x_{j+1},\ldots,x_{n-1},a_j] = 0 \quad\text{for $1 \leq j \leq k$ and all $x_l \in L$}.
\]
Then
\begin{equation}
\label{lk}
L_k = (L,\bracket_k,(\alpha_{k+1},\ldots,\alpha_{n-1}))
\end{equation}
is an $(n-k)$-ary Hom-Nambu algebra, where
\[
[x_{k+1},\ldots,x_n]_k = [a_1,\ldots,a_k,x_{k+1},\ldots,x_n]
\]
for all $x_l \in L$.  Moreover, if $L$ is multiplicative, then so is $L_k$.
\end{corollary}

\begin{corollary}
\label{cor3:lower}
Let $(L,\bracket,\alpha=(\alpha_1,\ldots,\alpha_{n-1}))$ be an $n$-ary Hom-Nambu-Lie algebra with $n \geq 3$.  Suppose for some $k \in \{1,\ldots,n-2\}$ there exist $a_i \in L$ for $1 \leq i \leq k$ satisfying
\[
\alpha_i(a_i) = a_i \quad \text{for $1 \leq i \leq k$}.
\]
Then $L_k$ in \eqref{lk} is an $(n-k)$-ary Hom-Nambu-Lie algebra.  Moreover, if $L$ is multiplicative, then so is $L_k$.
\end{corollary}

In Corollary \ref{cor3:lower} the assumption $\alpha_i(a_i) = a_i$ is automatically satisfied if $\alpha_i$ is the identity map.  In particular, it holds if $L$ is an $n$-ary Nambu-Lie algebra.  Therefore, we have the following special case of Corollary \ref{cor3:lower}.

\begin{corollary}
\label{cor4:lower}
Let $(L,\bracket)$ be an $n$-ary Nambu-Lie algebra with $n \geq 3$.  Let $a_1,\ldots,a_k \in L$ be arbitrary elements for some $k \in \{1,\ldots,n-2\}$.  Then $L_k = (L,\bracket_k)$ is an $(n-k)$-ary Nambu-Lie algebra, where
\[
[x_{k+1},\ldots,x_n]_k = [a_1,\ldots,a_k,x_{k+1},\ldots,x_n]
\]
for all $x_l \in L$.
\end{corollary}

The results in this section can be applied to any of the $n$-ary Hom-Nambu(-Lie) algebras from the previous sections.  Let us now discuss how Corollary \ref{cor1:lower} (with $n = 3$) can be used with a result in \cite{ams2} to create a different Hom-Lie algebra from a given one in a non-trivial way.

Let $(L,[,],\alpha)$ be a Hom-Lie algebra.  A \emph{trace function} is a linear form $\tau \colon L \to \bk$ such that $\tau([x,y]) = 0$ for all $x,y \in L$.  Given $L$ and $\tau$, suppose $\beta \colon L \to L$ is another linear map such that
\begin{equation}
\label{trace}
\begin{split}
\tau(\alpha(x))\tau(y) &= \tau(x)\tau(\alpha(y)),\\
\tau(\beta(x))\tau(y) &= \tau(x)\tau(\beta(y)),\\
\tau(\alpha(x))\beta(y) &= \tau(\beta(x))\alpha(y)
\end{split}
\end{equation}
for all $x,y \in L$.  Now consider the triple product
\begin{equation}
\label{taubracket}
[xyz]_\tau = \tau(x)[y,z] + \tau(y)[z,x] + \tau(z)[x,y],
\end{equation}
which was first introduced in \cite{almy}.  Under these hypotheses, it is proved in \cite{ams2} (Theorem 3.3) that
\begin{equation}
\label{ltau}
L_\tau = (L,[,,]_\tau,(\alpha,\beta))
\end{equation}
is a ternary Hom-Nambu-Lie algebra.  Using this construction of a ternary Hom-Nambu-Lie algebra from a Hom-Lie algebra in conjunction with Corollary \ref{cor1:lower} (with $n = 3$), we obtain the following result.  It gives a recipe for replacing the twisting map and the bracket in a Hom-Lie algebra in a non-trivial way.

\begin{corollary}
\label{cor5:lower}
Let $(L,[,],\alpha)$ be a Hom-Lie algebra, $\tau \colon L \to \bk$ be a trace function, and $\beta \colon L \to L$ be a linear map such that \eqref{trace} is satisfied.  If $a \in L$ satisfies $\alpha(a) = a$, then
\[
L_\tau' = (L,[,]_\tau',\beta)
\]
is a Hom-Lie algebra, where
\[
[x,y]_\tau' = \tau(a)[x,y] + [a,\tau(y)x - \tau(x)y]
\]
for $x,y \in L$.
\end{corollary}

\begin{proof}
By \cite{ams2} (Theorem 3.3) $L_\tau$ in \eqref{ltau} is a ternary Hom-Nambu-Lie algebra.  By Corollary \ref{cor1:lower} (with $n = 3$) it follows that $(L,[,]_\tau',\beta)$ is a binary Hom-Nambu-Lie algebra (i.e., a Hom-Lie algebra), where the reduced bracket $[,]_\tau'$ is given by
\[
\begin{split}
[x,y]_\tau'
&= [axy]_\tau\\
&= \tau(a)[x,y] + \tau(x)[y,a] + \tau(y)[a,x]\\
&= \tau(a)[x,y] + [a, \tau(y)x - \tau(x)y].
\end{split}
\]
In the last equality above, we used the anti-symmetry of the Hom-Lie bracket $[,]_\tau'$.
\end{proof}

%%==============%%
%%              %%
%%  References  %%
%%              %%
%%==============%%

\end{document}